\documentclass{article}
\usepackage{arxiv}
\usepackage{float}
\usepackage{enumitem}
\usepackage{enumerate}
\usepackage{tabularx}
\usepackage{tfrupee}
\usepackage{multicol}
\usepackage{amsmath}
\usepackage[utf8]{inputenc} 
\usepackage[T1]{fontenc}    
\usepackage{hyperref}       
\usepackage{url}            
\usepackage{booktabs}       
\usepackage{amsfonts}       
\usepackage{nicefrac}       
\usepackage{microtype}      
\usepackage{lipsum}
\usepackage{graphicx}
\graphicspath{ {./images/} }
\usepackage{array}
\newtheorem{thm}{Theorem}[section]
\newtheorem{defn}[thm]{Definition}

\title{Compositions of Knots Using Alexander Polynomial}

\date{}
\author {
 G INFANT GABRIEL \\
  PG and Research Department of Mathematics,\\
  Sri Ramakrishna College of Arts \& Science, \\
  Coimbatore-641006 ,\\
  Tamil Nadu, India.\\
\texttt{infantgabriel@srcas.ac.in} \\
   \And
Dr N UMA \\
  PG and Research Department of Mathematics,\\
  Sri Ramakrishna College of Arts \& Science, \\
  Coimbatore-641006 ,\\
  Tamil Nadu, India.\\
\texttt{uma.n@srcas.ac.in} \\
}

\begin{document}
\maketitle
\begin{abstract}
Knot theory is the Mathematical study of knots. In this paper we have studied the Composition of two knots. Knot theory belongs to Mathematical field of Topology, where the topological concepts such as topological spaces, homeomorphisms, and homology are considered. We have studied the basics of knot theory, with special focus on Composition of knots, and knot determinants using Alexander Polynomials. And we have introduced the techniques to generalize the solution of composition of knots to present how knot determinants behave when we compose two knots.	 

\end{abstract}


\section{Introduction}
Knot theory is a new kind of applicable Mathematics. \textbf {Edward E. David, Jr. \cite{David}} referred that utility of Mathematics conceived symbolically. The history of Knot theory starts from 19th century physics, with the work of Gauss on computing linking numbers in a system of linked circular wires.        \textbf{ D. Silver \cite {Silver}} also studied the knots and coined the word topology. In 1867, Kelvin’s vortex model of atom by \textbf {W.T Thompson \cite{Thompson} } was presented. The American Mathematician \textbf {J.W Alexander \cite{Alexandar}} was the first to suggest that knot theory is extremely important in the study of 3-dimensional topology which was further underlined by the German mathematician \textbf { H. Seifert.  Later K. Murasugi \cite{Murasugi}} studied the relationship between algebraic geometry and the knot theory. \\
Knot theory was first presented by the physicist and chemists \textbf {William Thomson} based on the study by \textbf { Baron Kelvin} who hypothesized that atoms of different elements can be defined by different knots. Through Thomson’s theory was later proved incorrect, his work inspired \textbf {Peter Trait}, who developed many concepts that are used today in application of knots theory to biology, chemistry and physics. \textbf {Colin Adams \cite{Colin}} suggested that the knot theory was used in modelling of DNA, the effects of enzymes and in statistical mechanics while examining the interaction between particles in a system. By using more theoretical models, Scientists and Mathematicians can make these concrete concepts to manipulate and work within.
And in this paper, we have studied knots that are embedded in $S^3$ . A knot can be projected onto a plane (or simply drew it on a paper). These projections are called knots diagrams (Figure 1). \\\\
To avoid ambiguity the following restrictions are considered
\begin{itemize}
	\item   At each crossing, the string segment passes over respectively under the other (this is usually done by drawing a gap in the bottom segment).
	\item 	Each crossing involves exactly two segments of the string.
	\item	The segments must cross transversely.
	\item 	At each crossing there is always one over strand and one understand. 
	\item 	An arc is a piece of the knot that passes from one undercrossing to another with only overcrossing in between (an unbroken line).
\end{itemize} 
The following diagrams are projections of the knots.
\begin{figure}[h]
			\centering
			\includegraphics[width=2in, height=2in]{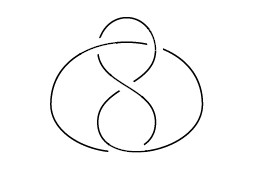}
			\caption{Figure 8 Knot}
		\end{figure}\\
		\begin{figure}[h]
		\centering
		\includegraphics[width=3in, height=2in]{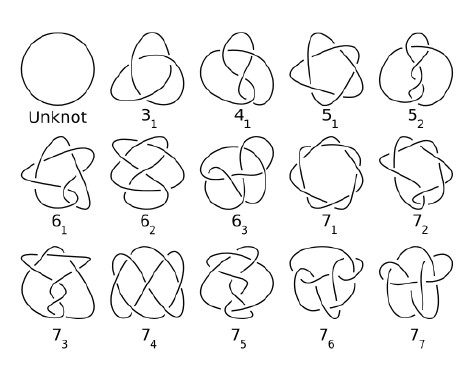}
		\caption{Knot Diagram}
	\end{figure}\\
\section{PRELIMINARIES}
Basic definitions and concepts are presented in this chapter
\begin{defn} \cite{Richard}
	\\
	K is a knot for the embedding $ \textit h:S^1\rightarrow  S^3 $  whose image is K.  If $K\subseteq S^3$ , then it is homeomorphic to the unit circle $S^1$  . The actual knot is a smooth embedding of the unit circle in $S^3$  . A knot can only consist of one component; a link on the other hand is a finite union of disjoint knots.
	\end{defn}
\begin{defn} \cite{Richard}
	\\
	Knots $K_1$  and $K_2$  in $S^3$  are \textbf equivalent if there exist a homeomorphism  $ \textit h:S^3\rightarrow S^3 $ such that $\textit h(K_1)=h(K_2)$
		\end{defn}
\begin{defn} \cite{Richard}
	\\
	If H is an isotopy between ambient spaces $\textit H:S^3\times I \times \rightarrow S^3$  , then H is an ambient isotopy. The knot can be deformed to any expected manner. The arcs can be bent and moved through space without passing through one another (knot can be shrunk or grown). Also, it is not permitted to pull the knot so tight that it unknots itself by disappearing into a point.
	\end{defn}
	\begin{defn} \cite{Mari}
\begin{itemize} 
	\item  \textbf {R1} allows to remove (or introduce) a twist in a diagram. The result is that the knot will have one fewer (or one more) crossing. 
	\item 	\textbf	{R2} lets separate two strings that lie on top of each other or vice versa. This will add or remove two crossings. 
	\item	\textbf {R3} allows moving a strand from one side of a crossing to the other. This also works if the strand is moved above the other two strands. 
\end{itemize} 
\end{defn}
R3 does not affect the number of crossings in the current projection. If a deformation of the diagram uses R2 and R3, referred as Regular Isotopy (or Planar Isotopy).The regular isotopy is an equivalence relation for the knot diagrams and is not defined for the knot embedding. The Reidemeister moves describes the procedures performed on diagrams
	\begin{figure}[h]
	\centering
	\includegraphics[width=4in, height=2in]{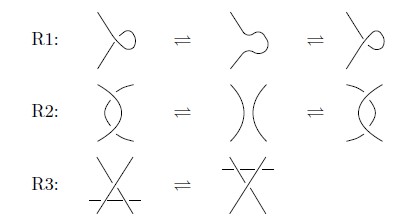}
	\caption{Reidemeister Moves Diagram}
\end{figure}
\begin{defn} \cite{Mari}
\\
	The crossing number, c(K) of a knot K,  is the \textbf {minimum number of crossings,} that occur in any diagram of K.
\end{defn}
\begin{defn} \cite{Mari}
\\
	A knot is \textbf oriented if it has an orientation assigned and refered by arrows.
	If a knot K has a given orientation, -K refers as \textbf {reverse orientation.} \\
	A knot is called invertible if K and -K are equivalent. All knots in Figure 2 are invertible.
		\begin{figure}[h]
		\centering
		\includegraphics[width=4in, height=1.5in]{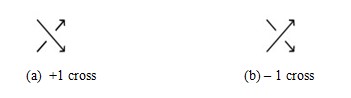}
		\caption{Positive and Negative Crossing}
	\end{figure}
	\end{defn}
\begin{defn} \cite{Mari}
\\
	The \textbf {writhe} of an oriented knot, w(K), is the sum of the crossings with signs as shown in Figure 4.
\end{defn}
\begin{defn} \cite{Mari}
\\
	A knot K is called \textbf {alternating} if its diagram, the undercrossing and overcrossings alternate around K.
\end{defn}
\begin{defn} \cite{Mari}
\\
	The \textbf {knot sum or connected sum $K_1 \# K_2$ }  is formed by placing two knots side by side, removing a small arc from each knot and then joining the knots together with two new arcs.
\end{defn}
\begin{defn} \cite{Mari}
\\
	A knot is called \textbf {composite} if it can be written as the sum of $K_1$ and $K_2$  , neither of which is the unknot and prime if it is not composite. Knot tables only show prime knots, see         Figure 2.
\end{defn}
\section{COMPOSITION OF KNOTS}
Given two projections of knots, a new knot obtained by removing a small arc from each knot projection and then connecting the four endpoints by two new arcs as in Figure 5. The resulting knot is the \textbf {composition of the two knots}. If we denote the two knots by the symbols  $K_{4_1}$ and $K_{3_1}$ , then their composition is denoted by $K_{4_1}$ \# $K_{3_1}$ . 
\begin{figure}[h]
	\centering
	\includegraphics[width=4in, height=2in]{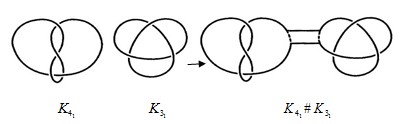}
	\caption{$K_{4_1}$ \# $K_{3_1}$ of two knots  $K_{4_1}$ and $K_{3_1}$ }
\end{figure}
The two projections should not overlap, and the two arcs are chosen to remove the outside of each projection to avoid any crossings i.e., they do not cross either the original knot projections or each other (Figure 6).
\begin{figure}[h]
	\centering
	\includegraphics[width=4in, height=2in]{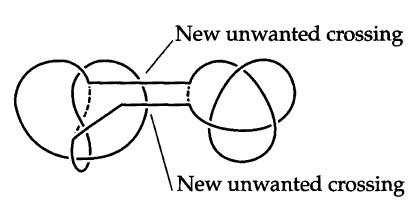}
	\caption{Not the Composition of $K_{4_1}$ and $K_{3_1}$ }
\end{figure}
\begin{defn}
	A knot a composite knot if it can be expressed as the composition of two knots, neither of which is the \textbf { trivial knot}. The knots that makeup the composite knot is called \textbf { factor knots.} \\
	The composition of a knot K with the unknot, the result is again K. (Figure 7). If a knot is not the composition of any two nontrivial knots, we call it a \textbf {prime knot.}
	\begin{figure}[h!]
		\centering
		\includegraphics[width=2in, height=2in]{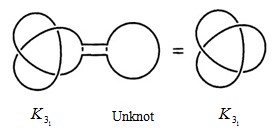}
		\caption{$K_{3_1}$ \#  Unknot gives $K_{3_1}$ }
	\end{figure}
	\end{defn} 
\begin{defn} 
	An \textbf {orientation} is defined by choosing a direction to travel around the knot. The orientation on $K_{3_1}$ matches the orientation on $K_{5_2} $ in $K_{3_1}\#K_{5_2}$   , resulting in an orientation for $K_{3_1}\#K_{5_2}$   , or the orientation on $K_{3_1}$ and $K_{5_2} $  do not match up in $K_{3_1}\#K_{5_2}$    .If the orientations do match up in all the compositions of the two knots then it will yield the same composite knot. If the orientations do match up in all the compositions of the two knots then it will yield the single composite knot.
		\begin{figure}[h]
		\centering
		\includegraphics[width=5in, height=1.5in]{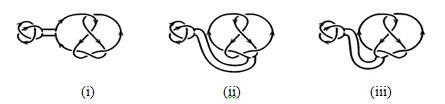}
		\caption{ (i) \& (ii)$K_{3_1}\#K_{5_2}$ Orientations match (iii) $K_{3_1}\#K_{5_2}$ Orientations match differ. }
	\end{figure}
\end{defn}
\section{THE ALEXANDER POLYNOMIAL OF COMPOSITION}
In 1928, \textbf {J.W Alexander \cite{Alexandar}} introduced polynomial invariant to compute the knot diagram entries of a matrix determinant called the Alexander polynomial. Alexander determined the entries of a matrix from the crossing and arcs of the diagram. The Alexander polynomial does not depend on the indexing of crossing and arcs. It does not also depend on the row and column eliminated from the crossing/arc matrix. The Skein relation discovered by John Conway in 1969 is to compute the Alexander polynomial. The Alexander polynomial is invariant up to multiplication by $\pm t^N$ where N is some integer. \\ \\
\textbf {Numerical Example:} \\
Composition of two knots  $K_{3_1}\#K_{1_1}$
\begin{figure}[h]
	\centering
	\includegraphics[width=4in, height=2in]{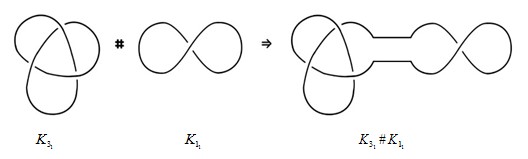}
\end{figure}
\\ \\ \\ \\ \\ Applying the condition of Alexander polynomial in above composition knot $K_{3_1}\#K_{1_1}$  we get,
\begin{figure}[h]
	\centering
	\includegraphics[width=4in, height=3in]{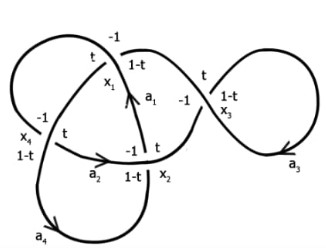}
\end{figure}
\begin{equation*}
	M_{{K_{3_1}}\#{K_{1_1}}}=\begin{pmatrix}
		1-t&0&-1&t\\
		t&1-t&0&-1\\
		0&-1&t-t^2&0\\
		-1&t&0&1-t\\
		\end{pmatrix}
\end{equation*}
The Alexander matrix $A_k$ is defined as the matrix \textit M by deleting row n and column n, where
\begin{equation*}
	A_{{K_{3_1}}\#{K_{1_1}}}=\begin{pmatrix}
		1-t&0&-1\\
		t&1-t&0\\
		0&-1&t-t^2\\
		\end{pmatrix}
\end{equation*}
The Alexander polynomial $ \Delta_k(t)$ of a Knot K is the determinant of its Alexander matrix, which is $\Delta_{{K_{3_1}}\#{K_{1_1}}} =2t-3t^2+3t^3-t^4$ \\
Finally, the resulting determinant is multiplied by $t^{-1}$ to normalize the Alexander polynomial in order to have a positive constant term. \\
Thus $2-3t+3t^2-t^3$ is the Alexander Polynomial of this Composition of Knot $K_{3_1}\# K_{1_1}$ .\\
We observe that Alexander polynomial for compositions of $K_{3_1}$ and $K_{1_1}$  is $2-3t+3t^2-t^3$.\\

Similarly the composition of$K_{3_1}$  and $K_{2_1} $ is $3-t^{-1}-5t+5t^2-t^3$  .
\begin{itemize}
	\item $K_{4_1} \rightarrow K_{1_1}, K_{2_1}$ are found to be $1-3t+4t^2-2t^3$   and  $2+t-8t^2+8t^3-4t^4+t^5$   respectively.
	\item $K_{5_1} \rightarrow K_{1_1}, K_{2_1}$ together to form $2-5t+2t^2-7t^3+5t^4-t^5$ and $3-t^{-1}-9t+12t^2-14t^3+10t^4+2t^5$.
	\item $K_{5_2} \rightarrow K_{1_1}, K_{2_1}$ are joined to be $2-t^{-1}-5t+7t^2-6t^3+2t^4$ and $1-t^{-1}-2t+7t^2-10t^3+6t^4-2t^5$.
	\item $K_{6_1}, K_{6_2},K_{6_3} \rightarrow K_{1_1}$ are found to be  $4-10t+4t^2+7t^3-8t^4 ,  2-2t-3t^2+3t^3+2t^4-2t^5$ and  $1-3t+2t^2+2t^3-2t^4$.
\end{itemize}
\begin{center}
    
	\begin{tabular}{|c|c|c|} \hline
	\textbf {S.No} & \textbf {Composition of two Knot} & \textbf {Composition Solution using AP} \\ \hline
	1 & $K_{3_1}\#K_{1_1}$ & $2-3t+3t^2-t^3$ \\ \hline
	2 &   $K_{3_1}\#K_{2_1}$ & $3-t^{-1}-5t+5t^2-t^3$ \\ \hline
	3 & $K_{4_1}\#K_{1_1}$ & $1-3t+4t^2-2t^3$ \\ \hline
	4 &  $K_{4_1}\#K_{2_1}$ & $2+t-8t^2+8t^3-4t^4+t^5$ \\ \hline
	5 & $K_{5_1}\#K_{1_1}$ & $2-5t+2t^2-7t^3+5t^4-t^5$ \\ \hline
	6 & $K_{5_1}\#K_{2_1}$ & $3-t^{-1}-9t+12t^2-14t^3+10t^4+2t^5$ \\ \hline
	7 & $K_{5_2}\#K_{1_1}$ &  $2-t^{-1}-5t+7t^2-6t^3+2t^4$ \\ \hline
	8 & $K_{5_2}\#K_{2_1}$ & $1-t^{-1}-2t+7t^2-10t^3+6t^4-2t^5$ \\ \hline
	9 & $K_{6_1}\#K_{1_1}$ &  $4-10t+4t^2+7t^3-8t^4$ \\ \hline
	10 &  $K_{6_2}\#K_{1_1}$ & $2-2t-3t^2+3t^3+2t^4-2t^5$ \\ \hline
	11 & $K_{6_3}\#K_{1_1}$ & $1-3t+2t^2+2t^3-2t^4$ \\ \hline
	\end{tabular}
\end{center}
 
\section{CONCLUSION}
Knots have always been an integral part of real life and have found Mathematical uses recently. As we can also visualize the application of compositions of knots in DNA replication, in this paper we have explained a brief about the composition of knots. Algebraic area is more important in knot theory and it involves drawing comparison between two knot compositions, hence we have applied the Alexander polynomial to get different polynomial for compositions of knots for $K_{3_1}\#K_{1_1},K_{3_1}\#K_{2_1},....etc., $ and in future we will extend the compositions of knots for $K_{6_2}\#K_{2_1},K_{7_1}\#K_{1_1},....etc., $


\begin{thebibliography}{1}
	
	\bibitem{Alexandar}J.W Alexandar, Topological invariants of knots and links, Transactions of the American Mathematical Society 30 (1928), 275-306.
	\bibitem{Colin} Colin Adams, “An Elementary Introduction to the Theory of Knots”, New York, WH Freeman, 1994.
	\bibitem{David}David, E.E. Jr. Renewing US mathematics: An agenda to begin the second century. Notices of the A.M.S. 35 (1988), 1119-1123.
	\bibitem{Mari}	Mari Ahlquist, “On Knots and DNA”, Department of Mathematics, Linköping University, LiTH-MAT-EX–2017/17–SE, December 2017.
	\bibitem{Murasugi}K. Murasugi, knot theory and its application, Birkhuser Boston, (1993).
	\bibitem{Richard}Richard H. Cromwell and Ralph H. Fox, Introduction to knot theory, Springer-Verlog; (1963).
	\bibitem{Silver}D. Silver, Knot theorys odd origins. American Scientist 94 (2006), 158-165.
	\bibitem{Thompson}W.T Thompson, Mathematical and Physical papers, III Cambridge U. Press (1890).
\end{thebibliography}
\end{document}